\documentclass[a4paper,12pt]{article}
\usepackage{authblk}
\usepackage{amsmath}
\usepackage{amsthm}
\usepackage{amssymb}
\usepackage[ruled,linesnumbered]{algorithm2e}
\usepackage{blindtext}
\usepackage{xcolor}
\usepackage{graphicx}
\usepackage{caption}
\DeclareCaptionLabelFormat{algnonumber}{Algorithm}

\usepackage{wrapfig}

\usepackage{multicol}
\usepackage{units}
\usepackage{ifpdf}
\usepackage{cite}

\usepackage{paralist}
\usepackage{inputenc}
\usepackage[left=1in,top=1in,right=1in,bottom=1in]{geometry}
\usepackage{indentfirst}
\title{\textbf{The 2-center Problem in Maximal Outerplanar Graph} }
\author{Hsiu-Fu,Yeh}

\affil{Department of Computer Science and Information Engineering\\
National Chung Cheng University, Chiayi, Taiwan
a007123456@gmail.com}
\date{}
\begin{document}

\maketitle
%
%
\begin{abstract}
    We consider the problem of computing 2-center in maximal outerplanar graph. In this problem, we want to find an optimal solution where two centers cover all the vertices with the smallest radius. 

    We provide the following result. We can compute the optimal centers and the optimal radius in $O(n^2)$ time for a given maximal outerplanar graph with $n$ vertices. We try to let the maximal outerplanar graph be cut into two subgraphs with an internal edge, each center will cover vertices which are continuous.
\end{abstract}

%
%

\section{Introduction}

An outerplanar grpah is first characterized by Chartrand, Gary and Harary, Frank\cite{Permutation}, and then Maciej M.Sysło\cite{Characterizations} defined some characterizations of the outerplanar graph. An outerplanar graph is a planar segmentation of a polygon. The center of a graph is defined to be the subgraph induced by the set of vertices that have minimum eccentricities. (i.e., minimum distance to the most distant vertices). Various types of centers have been defined for graphs, involving extreme values of distance, weighted, etc., defined for vertices of a graph. We will consider the most common of these, that based on distances between vertices. In this paper, we always assume the graph is simple which is without loops and multiple edges.

The maximal outerplanar graph is an outerplanar graph that can not have any additional edges added to it while preserving outerplanarity. Every bounded face of a maximal outerplanar graph is a triangle. In this paper, we consider the 2-center problem in a maximal outerplanar graph that the centers set has two centers with the smallest radius $r$ to cover all vertices.

The center problem is want to find an optimal center with the smallest radius $r$ for the eccentricities and the center can cover all vertices. Andrzej Proskurowski \cite{Centers} found the possibles for the centers of an outerplanar graph. In previous works on maximal outerplanar graphs, Arthur M.Farley and Andrzej Proskurowski \cite{Computation} defined the edge eccentricities for the edges of an outerplanar graph. They presented an algorithm that efficiently computes these edge eccentricities. And provided a linear time algorithm computation of the center.

The multiple-centers problem is provided by Minieka, Edward\cite{m-center}. A k-center set of a graph is any set of k points, that minimizes the maximum distance from a vertex to its nearest center; The description of the problem was defined by Khuller, Samir and Sussmann, Yoram J.\cite{Capacitated}.Given a graph $G$ = ($V$, $E$), find a subset $S \subseteq V$ of size $k$ such that each vertex in $V$ is close to some vertex in $S$, the objective function is defined as follows: 
\[
\mathop{min}\limits_{S \subseteq V} \mathop{max}\limits_{u \in V} \mathop{min}\limits_{v \in S} d(u,v)
\]

\section*{\large{The result}}

We focus on the centers of maximal outerplanar graph, it will cut the vertices into two continuous sets with it’s order. There exists an internal edge to cut the vertices into two continuous sets. By the way, we combine the algorithm that can compute the radius of an outerplanar graph was provided by Arthur M.Farley and Andrzej Proskurowski \cite{Computation}. We provide an algorithm to calculate the optimal radius and locate the centers of a maximal outerplanar graph for 2-center problem in $O(n^2)$ time . We will begin with the following definitions.

%
%
\section{Notation and Definition}
The planar graph is outerplanar if and only if it can be embedded in the plane with all vertices lie on the same face. Usually assume this face is the exterior. An edge cuts the graph into two sides, which has one side is empty (or called exterior side) are called the bounded edge. Others are internal edges. 

\paragraph{Definition 1. Distance}
  
    In a connected graph $G$, the path $p(v,w)$ is the shortest path from the vertex $v$ to vertex $w$. The path $p(v,S)$ is the farthest path in every $p(v,w)$ for all vertices $w$ in set $S$. The distance, $d(v,w)$ , is a value which between vertices $v$ and $w$ is the length of the shortest path joining them, or say the absolute value $|p(v,w)|$ . And $d(v,S)$ is the maximal value of between vertex $v$ and other vertices in the set $S$. For any vertex $v$,$x$,$w$ in $G$, $d(v,w) \leq \ d(v,x) + d(x,w)$.
    
\vspace{1ex}    
On a maximal outerplanar graph (called mop before), the ends of an edge can always find a vertex on non-exterior side such that the vertex is a neighbor of the ends. For an edge e in a mop will split the graph into two mops of the sides of $e$. There is no edge between the vertex belongs to two mini-mops without $e$.
    
\paragraph{Definition 2. Center and Radius}
  
    In a connected graph $G$, the eccentricity, $e(v)$ , of a vertex $v$ is the largest distance from v to any vertex of G; The radius, $r(G)$ or $r$ in simple, is the minimum eccentricity of vertices in $G$; The centers, $C$, are the set of the vertices, which have minimum eccentricity.
 
\paragraph{Definition 3. Center's cover region}
    For a center $c$ of a connected graph $G$ with radius $r$, the cover region of $c$ is defined that $R_c$ induced all the vertices v if $d(c,v) \leq  r$.
    
\centerline{\includegraphics{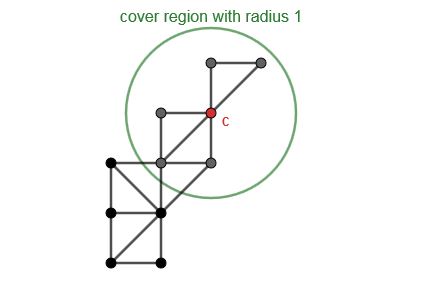}}

\centerline{\label{fig:cover_region}Fig.1 A cover region with a center $c$.}

\paragraph{Definition 4. Center's cover set}    
     Let the vertices be covered by center $c$ and call them in the $c$ cover set. A vertex will only be in one center's cover set whether it can be covered by another center with the same radius or not. It can change into another cover set if it can be covered by that center.

\vspace{1em}  
In this paper we need two centers to cover the vertices. We want to differentiate the vertex that is covered by which center. Use the cover sets can check the relation of centers and vertices efficiently.
%
%
\section{Continuous cover set}

We consider the external subclass of outerplanar graphs known as maximal outerplanar graphs, or “mops” in this paper.

\paragraph{Definition 5. Continuous cover set}    
    For a Hamilton cycle of a mop G, labeling the vertices with their order on it. If the orders of vertices in a cover set is continuous, called that cover set is continuous.

In this part, we show how to put the vertices into different center's cover set without increasing the radius if the centers are optimal.

\paragraph{Lemma 1.} 
The cover set of an optimal 2-center covered has continuous vertex order with the Hamilton cycle for a mop.

\vspace{1em}
\noindent \textbf{$\proofname.$}
Assume a mop $G$ and the two centers of cover sets are noncontinuous with an optimal radius $r$. Let the cover set called $S_1$, $S_2$, $S_3$ and $S_4$, let $S_1$, $S_3$ is covered by a center $c_1$ and $S_2$, $S_4$ is covered by another center $c_2$. Let $c_1$ in $S_1$ and $c_2$ in $S_2$. As $G$ is a mop, there is only one edge connected between $S_1$ and $S_3$ or $S_2$ and $S_4$. Without loss of generality, assume $S_1$ and $S_3$ is connected with edge $e_{con}$.

\centerline{\includegraphics{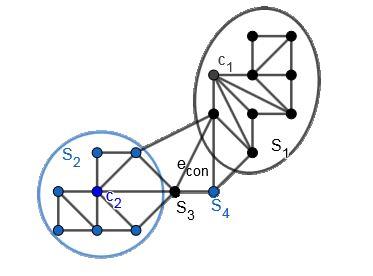}}

\centerline{\label{fig:cover_set}Fig.2 A sample of cover set $S_1$ to $S_4$ with radius 2.}
\vspace{1em}
Consider path $p(c_2,S_4)$ first, the value of $d(c_2,S_4)$ is less than $r$.And take a vertex $v_x$ suitable, let the path $p(c_2,S_4)$ be equal to $p(c_2,v_x)$ + $p(v_x,v_y)$ + $p(v_y,S_4)$, which $v_x$ belongs to $S_2$ , $v_y$ is one of the ends of $e_{con}$, and $v_x$ has a neighbor is $v_y$.

Without loss of generality, let $v_y$ belong to $S_1$ and another end of $e_{con}$ called $v_{y'}$. We will consider the path and distance between $c_2$ and $S_3$, $p(c_2,S_3)$ and $d(c_2,S_3)$. We know that, $d(c_2, S_2) \leq \  r $ and $ d(c_2, S_4) \leq \  r  $

$S_2$ and $S_4$ are separated by the edge $e_{con}$. The shortest path between $c_2$ and every vertex in $S_4$ must pass through one of vertices with its end.

\[
p(c_2, S_4) = p(c_2, v_x) + p(v_x, v_{y}) + p(v_{y}, S_4) 
\]
\subparagraph{(i)} If $d(c_2, S_3) > r$ 

\begin{equation}
 d(c_2, S_3) = d(c_2, v_x) + d(v_x, v_y) + d(v_y, v_{y'}) + d(v_{y'}, S_3) > r. 
\end{equation}

\begin{equation}
 d(c_1,S_3) = d(c_1,v_y) + d(v_y,v_{y'}) + d(v_{y'},S_3) \leq r.
\end{equation}

by (1) and (2), we know that $d(c_2,v_y) \geq  d(c_1,v_y)$.

\vspace{1em}
$
d(c_1,S_4) = d(c_1,v_y) + d(v_y,S_4) < d(c_2,v_y)+ d(v_y,S_4) = d(c_2,S_4) \leq r
$

We can change $S_4$ into $c_1$ cover set, and $S_1 \cup S_3 \cup S_4$ is a continuous cover set.

\subparagraph{(ii)}If $d(c_2,S_3) \leq r$

Just change the whole set $S_3$ into $c_2$ cover set, and then $S_2 \cup S_3 \cup S_4$ is a continuous cover set. Hold

\vspace{1em}
By \textbf{(i)} and \textbf{(ii)}, we know that there exist a least one separation of vertices of the optimal centers that is continuous order. If we want to find an optimal solution for the problem, we can cut the vertex with the order into two mops and then calculate the optimal center and radius of the subgraphs.

%
%

\section{Separation of vertices}

Although we can separate the vertex into two sets by the order, if want to calculate all the possible sets of vertices is inefficiency. So we try to use the edge in the graph to decrease the possible.

\paragraph{Lemma 2.}  A least one optimal center set exists that is separated by an internal edge. 

For a mop $G$ have $n$ vertices. Let the optimal center sets $S_1$, $S_2$ with centers $c_1$, $c_2$. $S_1 = \{ v_1,v_2...v_i\} $ and $S_2 = \{ v_{i+1},v_{i+2},...,v_n\} $. There exist an edge $e(v_1,v_i)$ or $e(v_{i+1}, v_n)$. And the radius is $r$.

\vspace{1em}
\noindent\textbf{$\proofname.$}
\textbf{Hypothesizing the mop $G$ have no edge connected $v_1$ and $v_i$ and also $v_{i+1}$ and $v_n$.}(i.e. $e(v_1,v_i)$ and $e(v_{i+1}, v_n)$ is not exist.)

Assume a mop $G$ have the optimal centers $c_1$ and $c_2$ with radius $r$.
\vspace{1ex}
First, let $c_1$ have a continuous cover set as large as possible. The cover set of $c_1$ is $v_1,v_2...v_i$ that is labeled clockwise on bounded edge.
$d(c_1,v_k)\leq r$, for $1\leq k \leq i$. Maybe some vertices can be covered by $c_1$ but noncontinuous, put them with other vertices in $c_2$ cover set by lemma 1.

\subparagraph{case 1 :}Both $v_1$ and $v_i$ are not in $c_2$ cover region.
\[
d(c_2,v_1)>r\;\; {\rm and}\;\; d(c_2,v_i)>r
\]

We know $v_n$ and $v_{i+1}$ are in $c_2$ cover set. So the shortest path between $c_2$ and $v_1$ is $p(c_2,v_1) = p(c_2,v_n) + p(v_n,v_1)$. Also $p(c_2,v_i) = p(c_2,v_{i+1})+ p(v_{i+1},v_i)$.
s.t. $d(c_2,v_n)$ and $d(c_2,v_{i+1})$ is equal to $r$.

$c_2$ and $v_1$ are separated by an edge which one of the ends is $v_n$.
That means exist an edge between $v_n$ and a vertex $v_b$ separates $c_2$ and $v_1$.
There also exist an edge form $v_{i+1}$ separates $c_2$ and $v_i$. Called the end of the edge $v_c$.
One of $v_b$ or $v_c$ belongs to $c_1$ cover set and the other belongs to $c_2$ cover set.

Without loss of generality, let $v_b$ in $c_2$ cover set. If $v_b$ in $c_2$ cover set and not $v_{i+1}$, then the path of $c_2$ to $v_{i+1}$ will rewrite into
\[
d(c_2,v_{i+1}) = r = d(c_2,v_b) + d(v_b,v_{i+1}) 
\]
and
\[
d(c_2,v_n) = d(c_2,v_b) + d(v_b,v_n) = d(c_2,v_b) + 1 = r
\]
\[
s.t\;\; v_b = v_{i+2}
\]
If an edge connects $v_n$ and $v_{i+2}$ means $v_{i+1}$ does not connect any vertex in $c_2$ cover set without $v_{i+2}$ or $v_n$. So $v_c$ belongs to $c_1$ cover set. Hold

Then, consider the sides of the edge $e(v_{i+1},v_c)$. There is no edge between the two sides.
We know $c_1$ is lain on ${v_1,v_2,...,v_c}$. 
$p(c_1,v_i)$ will equal to $p(c_1,v_{i+1}) + p(v_{i+1},v_i)$ or $p(c_1,v_c) + p(v_c,v_i)$ , the first one is contradiction, trivially.
So take $p(c_1,v_i) = p(c_1,v_c) + p(v_c,v_i)$

\[r = d(c_1,v_c) + d(v_c,v_i) \geq d(c_1,v_c) + d(v_c,v_{i+1})   \\
 = d(c_1,v_c) + 1
\]


It will let $v_{i+1}$ can be covered by $c_1$ with radius $r$. It's also contradiction.


\subparagraph{case 2 :} There exist some vertices in $c_1$ cover set can be covered by $c_2$

There exist some vertices in the overlapping of two center cover regions. 
Let the continuous region $R_2$ collect vertex $v_{r_2}$ that cover by both $c_1$ and $c_2$ from $r_2 = 1,2,...,k$, if $d(c_1,v_{r_2}) \leq r, {\rm and}\  d(c_2,v_{r_2}) \leq r$, and the continuous region $R_4$ collect the vertices $v_{r_4}$ that cover by both $c_1$ and $c_2$ from $r_4 = i,i-1,...,i-j$, too.

We cut the other vertices into regions $R_1$ and $R_3$ which are covered only by $c_1$ or $c_2$. Next, will show that there exist an edge between $R_2$ and $R_4$. Think about the distance of centers and the region $R_2$,$R_4$. 
\[
d(c_1,R_2) = d(c_1,R_4) = r = d(c_2,R_2) = d(c_2,R_4)
\]

If there have no edge between $R_2$ and $R_4$ will exist an edge between $R_1$ and $R_3$ called $e(v_{R_1},v_{R_3})$ with two ends $v_{R_1},v_{R_3}$. The edge will separate $c_2$ and one set of $R_2$ or $R_4$. Without loss of generality, let $c_2$ and $R_2$ be separated. 

Then $p(c_2,R_2) = p(c_2,v_{R_1}) + p(v_{R_1},R_2)$ or $p(c_2,v_{R_3}) + p(v_{R_3},R_2)$ which is smaller.


 \vspace{1ex}

Whether the path is through $v_{R_1}$ or $v_{R_3}$. All the vertices in $R_2$ are in $c_2$ cover region. $p(c_2,v_{R_1}) \leq  p(c_2,v_{R_3})+ p(v_{R_3},v_{R_1})$.
So that $r = d(c_2,R_2) \geq d(c_2,v_{R_1})$. $v_{R_1}$ will be in $c_2$'s cover region. It's contradiction. Such that there exist an edge between $R_2$ and $R_4$. 
If an edge between $R_2$ and $R_4$ is exist. We can get the vertex of the cover sets and rename the label without increasing the radius.

 \vspace{1ex}
By \textbf{case (1) {\rm and} (2)}, there is a way to label the vertex so that one of $e(v_1,v_i)$ or $e(v_{i+1}, v_n)$ will exist.

\paragraph{Lemma 3.}  Exist an internal edge cut the mop $G$ into two mini-mops.

%
%

\section{An edge eccentricities algorithm of mops}
This algorithm was presented by Arthur M.Farley and Andrzej Proskurowski \cite{Computation} in 1980 for the computation of the center and diameter of outerplanar graphs.

\paragraph{Lemma 4.1}
Given an edge $p = (u, v)$ of a mop with a non-empty side $S$, let $e_1$, $e_2$
and $e_u$ represent the values of $e(b, w, S_2)$, $e(b, v, S_2)$ and the eccentricity of u in the graph $S_2 \cup \{ u, v, w\} $, respectively. Then, the value of $e_u$ is $-(1 + e_2)$ if $e_2 > 0$, and $|e_2|$ otherwise.\cite{Computation}


\centerline{\includegraphics{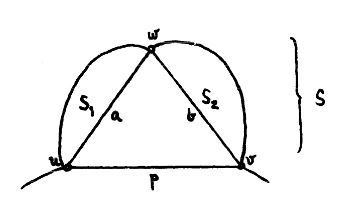}}
\centerline{\label{fig:egde_ecc}An edge $p = (u, v)$ with a nonempty side $S$.}

\paragraph{Lemma 4.2}
Given an edge $p = (u, v)$ with a non-empty side $S$, let $e_3$ and $d_1$
represent the values $e(a, u, S_1)$ and $e(p, u, S)$, respectively. Let $e_u$ be the eccentricity of $u$ in the graph $S2 \cup \{u, v, w\} $ as in $Lemma$ 3.1. Thus, the value of $d_1$ is $|e_3|$ if $|e_3| > e_u$. and  $e_u$ otherwise.\cite{Computation}

\paragraph{Lemma 4.3}
Given an edge $p = (u, v)$ with a non-empty side $S$, the eccentricity
$d_2 = e(p, v, S)$ is $|e(b, v, S_2)|$ if $|e(b, v, S_2)| \geq |q|$,and $|e(b, v, S_2)|$ otherwise,where $q$ is $-(1 + e(a, u, S_1) $if $e(a, u, S_1) > 0$, and $|e(a, u,S_1)|$ otherwise.\cite{Computation}

\vspace{1ex}
 As mentioned above, given a mop $G$, we can calculate the eccentricities of its edges as follows:

\textbf{(1)} For all edges $p = (u, v)$ on the Hamiltonian cycle of $G$, assign the value $-1$ to $e(p, x, \emptyset )$,where $x \in (u, v)$.

\textbf{(2)} For each triangle $(u, v, w)$ such that values $e(a, u, S_1)$, $e(a, w, S_1)$, $e(b, v, S_2)$, and $e(b, w, S_2)$ are defined, assign values of $e(p, u, S)$ and $e(p, v, S)$.
according to the rules specified in $Lemmas$ 4.1, 4.2 and 4.3.

\paragraph{Lemma 4.4} The eccentricity of a vertex $v$ is equal to the maximum absolute
value of the two pertinent eccentricities of any edge incident with $v$.\cite{Computation}

\paragraph{Algorithm 1} Given an outerplanar graph $G$ with known edge eccentricities,
compute the radius, and center of $G$ as follows\cite{Computation}:

(1) For each vertex $v$ of $G$ select an edge $p$ incident with $v$ and set eccentricity of $v$ to the maximum of the absolute values of the two pertinent edge eccentricities of $p$.

(2) Set the center of $G$ to be the set of vertices whose eccentricities equal the radius of $G$.

\vspace{1ex}
The time complexity of the algorithms for computing all vertex eccentricities in an outerplanar graph presented is $O(n)$. 

%
%

\section{Compute the radius and the 2-center on mops}
As mentioned previously, we can compute the optimal radius and the location of 2-center for a maximal outerplanar graph.
\paragraph{Algorithm 2} Given a mop $G$ with $n$ vertex, compute the radius, and 2-center of $G$ as follows:

\textbf{(1)} For each internal edge $e$ of $G$, takes the vertices induced one side of $e$ is a subgraph of $G$, is $G_{e_\alpha} $. And let another subgraph induces other vertices called $G_{e_\beta} $

\textbf{(2)} Calculate all the edge eccentricities in every $G_{e_\alpha}$ and $G_{e_\beta}$ use Algorithm 1, by Arthur M.Farley and Andrzej Proskurowski. Take the minimal eccentricities as radius in each subgraph. 

\textbf{(3)} Set the radius $r_e$ of edge $e$ to the maximal radius of $G_{e_\alpha}$ and $G_{e_\beta}$ 

\textbf{(4)} Set the optimal radius $r$ of $G$ to the minimal value of each $r_e$

\textbf{(5)} Set the 2-center of $G$ to be the centers, one in $G_{e_\alpha}$ and another in $G_{e_\beta}$ which the radius $r_e$ of edge $e$ equal to the optimal radius $r$.

The time complexity for computing of \textbf{(2)} is $O(n)$. The numbers of all internal edges are equal to $n-3$, every edge will be computed by \textbf{(2)} twice for each side of the edge, and \textbf{(3)}, \textbf{(4)}, \textbf{(5)} take linear time before all computed of \textbf{(1)} and \textbf{(2)}. So the time complexity of the algorithms is $O(n^2)$.

\renewcommand{\thealgocf}{2}

\begin{algorithm}

\caption{calculate optimal radius in maximal outerplanar graph}
\KwData{Maximal outplanar graph G, vertices set V, edges set E}
\KwResult{Optimal radius r, a centers set C}
\For{$e$ in $E$}{
$G_{e_\alpha} \leftarrow $ one side of $e$\;
$G_{e_\beta} \leftarrow $ the other side of $e$\;
\uIf{$G_{e_\alpha}$ and $G_{e_\beta}$ both is a mop}{
$r_\alpha \leftarrow $ the radius of $G_{e_\alpha}$; \qquad by algorithm 1\\    
$r_\beta \leftarrow $  the radius of $G_{e_\beta}$;  \qquad by algorithm 1\\  
$r_e \leftarrow $ max\{$r_\alpha$, $r_\beta$\} \;
}
}
$r \leftarrow \mathop{min}\limits_{e \in E} r_e$ if $r_e$ exist \;
\For {every $r_e$}{
\uIf{$r_e = r$}
{$C \leftarrow$ the center of $G_{e_\alpha}, G_{e_\beta}$\;
break\;
}
}

\end{algorithm}

\paragraph{Theorem 4} A 2-center problem with a given maximal outerplanar graph can be computed in $O(n^2)$ time.

%
%
\section{Conclusion}
In this paper, we provide the $O(n^2)$-algorithm to compute the optimal radius of a mop and find the 2-center. There are some extensions of this work in the future, e.g., a n-centers problem in a mop, or considering the edges and vertices with weighted. 

\bibliographystyle{unsrt}

\bibliography{ref}

\begin{thebibliography}{1}

\bibitem{Permutation}
Gary Chartrand and Frank Harary.
\newblock Planar {Permutation} {Graphs}.
\newblock {\em Annales de l'I.H.P. Probabilit\'es et statistiques},
  3(4):433--438, 1967.

\bibitem{Characterizations}
Maciej~M. Sysło.
\newblock Characterizations of outerplanar graphs.
\newblock {\em Discrete Mathematics}, 26(1):47--53, 1979.

\bibitem{Centers}
Andrzej Proskurowski.
\newblock Centers of maximal outerplanar graphs.
\newblock {\em Journal of Graph Theory}, 4(1):75--79, 1980.

\bibitem{Computation}
Arthur~M. Farley and Andrzej Proskurowski.
\newblock Computation of the center and diameter of outerplanar graphs.
\newblock {\em Discrete Applied Mathematics}, 2(3):185--191, 1980.

\bibitem{m-center}
Edward Minieka.
\newblock The m-center problem.
\newblock {\em SIAM Review}, 12(1):138--139, 1970.

\bibitem{Capacitated}
Samir Khuller and Yoram~J. Sussmann.
\newblock The capacitated k-center problem.
\newblock {\em SIAM Journal on Discrete Mathematics}, 13(3):403--418, 2000.

\end{thebibliography}

\end{document}